%\plain tex
\relax 
\global \edef \factsub {{0.1}}
\global \edef \factfewvars {{0.3}}
\global \edef \factdec {{0.4}}
\global \edef \sectminimal {{1}}
\global \edef \qsectminimal {{4}}
\global \edef \lmnosf {{1.1}}
\global \edef \propnosf {{1.2}}
\global \edef \propprimesnosf {{1.3}}
\global \edef \propprimes {{1.4}}
\global \edef \propminprime {{1.5}}
\global \edef \propmincomps {{1.6}}
\global \edef \propnotrad {{1.7}}
\global \edef \sectdblexpemb {{2}}
\global \edef \qsectdblexpemb {{9}}
\global \edef \qcalc {{10}}

%%%
\magnification=\magstep1
\baselineskip=13pt
\overfullrule=0pt

\font\footfont=cmr8
\font\footitfont=cmti8
\font\sc=cmcsc10 at 12pt
\font\large=cmbx12 at 14pt

%%%%%%%%%%%%%%%%%%% for math symbols not in plain tex %%%%%%%%%%%%%%%%%%%%
\font\tenmsb=msbm10
\font\sevenmsb=msbm7
\font\fivemsb=msbm5
\newfam\msbfam
\textfont\msbfam=\tenmsb
\scriptfont\msbfam=\sevenmsb
\scriptscriptfont\msbfam=\fivemsb
\def\hexnumber#1{\ifcase#1 0\or1\or2\or3\or4\or5\or6\or7\or8\or9\or
	A\or B\or C\or D\or E\or F\fi}

\mathchardef\subsetneq="2\hexnumber\msbfam28
%%%%%%%%%%%%%%%%%%%%%%%%%%%%%%%%%%%%%%%%%%%%%%%

%\def\emptyset{\{\}}

\def\hht{\hbox{ht}\hskip 0.3em}

\def\qed{\hbox{\quad \vrule width 1.6mm height 1.6mm depth 0mm}\vskip 1ex}
\def\eqed{\hbox{\quad \vrule width 1.6mm height 1.6mm depth 0mm}}
\def\mod{\hskip 0.3em\hbox{modulo}\hskip 0.3em}
\def\D{D}
\def\C{C}
\def\E{E}
\def\F{F}

\def\mvrule{\vrule height 3.2ex depth 1.6ex width .01em}

\newcount\tableco \tableco=0
\def\tabledo{\global\advance\tableco by 1{\the\tableco}}
\def\tabel#1{{\global\edef#1{\the\tableco}}}

\def\label#1{{\global\edef#1{\the\sectno.\the\thmno}}\ignorespaces}
\def\pagelabel#1{\immediate\write\isauxout{\noexpand\global\noexpand\edef\noexpand#1{{\the\pageno}}}{\global\edef#1{\the\pageno}}\ignorespaces}

\def\isnameuse#1{\csname #1\endcsname}

\def\issecond#1#2{#2}
\def\isifundefined#1#2#3{
	\expandafter\ifx\csname #1\endcsname\relax #2
	\else #3 \fi}
\def\pageref#1{\isifundefined{is#1}
	{{\bf ??}\message{Reference `#1' on page [\number\count0] undefined}}
	{\edef\istempa{\isnameuse{is#1}}\expandafter\issecond\istempa\relax}}

\def\today{\ifcase\month\or January\or February\or March\or
  April\or May\or June\or July\or August\or September\or
  October\or November\or December\fi
  \space\number\day, \number\year}

\newcount\sectno \sectno=0
\newcount\thmno \thmno=0
\def \section#1{\vskip 1.2truecm
	\global\advance\sectno by 1 \global\thmno=0
	\noindent{\bf \the\sectno. #1} \vskip 0.6truecm}
\def \thmline#1{\vskip 15pt
	\global\advance\thmno by 1
	\noindent{\bf #1\ \the\sectno.\the\thmno:}\ \ %
	%\noindent{\bf #1\ \the\thmno:}\ \ %
	\bgroup \advance\baselineskip by -1pt \it
	\abovedisplayskip =4pt
	\belowdisplayskip =3pt
	\parskip=0pt
	}

\def \thm{\thmline{Theorem}}

\def \endb{\egroup \vskip 1.4ex}

\def \lemma{\thmline{Lemma}}

\def \fact{\advance\thmno by 1\item{\bf\the\sectno.\the\thmno:}}
\def \prop{\thmline{Proposition}}
\def \proof{\smallskip\noindent {\sl Proof:\ \ }}

%%%%%%%%%%%%%%%%%%%%%%%%%%%%%%%%%%%%%%%%%%%%%%%%%%%%%%%%%
% Forward referencing

% The following will put labels also in a file so I can do forward referencing.
% I need an open file for that -- see below.
\def\label#1{\unskip\immediate\write\isauxout{\noexpand\global\noexpand\edef\noexpand#1{{\the\sectno.\the\thmno}}}
    {\global\edef#1{\the\sectno.\the\thmno}}\unskip\ignorespaces}
\def\sectlabel#1{\immediate\write\isauxout{\noexpand\global\noexpand\edef\noexpand#1{{\the\sectno}}}
    {\global\edef#1{\the\sectno}}}

\def\tabel#1{\hbox to 0pt{\hskip -6em\string#1}\unskip\immediate\write\isauxout{\noexpand\global\noexpand\edef\noexpand#1{{\the\tableco}}}\unskip}
%\def\tabel#1{\unskip\immediate\write\isauxout{\noexpand\global\noexpand\edef\noexpand#1{{\the\tableco}}}\unskip}

% Opening a file for forward referencing:
\newwrite\isauxout
\openin1\jobname.aux
\ifeof1\message{No file \jobname.aux}
       \else\closein1\relax\input\jobname.aux
       \fi
\immediate\openout\isauxout=\jobname.aux
\immediate\write\isauxout{\relax}
%%%%%%%%%%%%%%%%%%%%%%%%%%%%%%%%%%%%%%%%%%%%%%%%%%%%%%%%%

\ %
\vskip 3ex
\centerline{\large The minimal components of the Mayr-Meyer ideals}

\vskip 4ex
\centerline{\sc Irena Swanson}
\footnote{}{The author thanks the NSF for partial support
on grants DMS-0073140 and DMS-9970566.}
\centerline{\sc \today}
%\centerline{\sc \today, \number\time\ minutes past midnight}
\unskip\footnote{ }{{\footitfont 1991 Mathematics Subject Classification.}
13C13, 13P05}
\unskip\footnote{ }{{\footitfont Key words and phrases.}
\footfont Primary decomposition, Mayr-Meyer, membership problem, complexity of ideals.}
%}

\vskip 0.5cm

Grete Hermann proved in [H] that for any ideal $I$
in an $n$-dimensional polynomial ring over the field of rational numbers,
if $I$ is generated by polynomials $f_1, \ldots, f_k$ of degree at most $d$,
then it is possible to write $f = \sum r_i f_i$
such that each $r_i$ has degree at most $\deg f + (kd)^{(2^n)}$.
Mayr and Meyer in [MM] found (generators) of a family of ideals
for which a doubly exponential bound in $n$ is indeed achieved.
Bayer and Stillman [BS] showed that for these Mayr-Meyer ideals
any minimal generating set of syzygies
has elements of doubly exponential degree in $n$.
Koh [K] modified the original ideals to obtain homogeneous quadric ideals
with doubly exponential syzygies and ideal membership equations.

Bayer, Huneke and Stillman asked whether
the doubly exponential behavior is due to the number of minimal
and/or associated primes,
or to the nature of one of them?
This paper examines the minimal components and minimal primes of the Mayr-Meyer ideals.
In particular,
in Section~\sectdblexpemb\ it is proved that
the intersection of the minimal components of the Mayr-Meyer ideals
does not satisfy the doubly exponential property,
so that the doubly exponential behavior of the Mayr-Meyer ideals
must be due to the embedded primes.

The structure of the embedded primes of the Mayr-Meyer ideals is examined in [S2].

There exist algorithms for computing primary decompositions of ideals
(see Gianni-Trager-Zacharias [GTZ], Eisenbud-Huneke-Vasconcelos [EHV],
or Shimoyama-Yokoyama [SY]),
and they have been partially implemented
on the symbolic computer algebra programs Singular and Macaulay2.
However,
the Mayr-Meyer ideals have variable degree
and a variable number of variables over an arbitrary field,
and there are no algorithms to deal with this generality.
Thus any primary decomposition of the Mayr-Meyer ideals
has to be accomplished with traditional proof methods.
Small cases of the primary decomposition analysis
were partially verified on Macaulay2 and Singular,
and the emphasis here is on ``partially":
the computers quickly run out of memory.
 %(e.g., $d = 2$ and $n = 2, 3$)

The Mayr-Meyer ideals are binomial,
so by the results of Eisenbud-Sturmfels in [ES]
all the associated primes themselves are
also binomial ideals.
It turns out that many minimal primes are even monomial,
which simplifies many of the calculations.

The Mayr-Meyer ideals depend on two parameters,
$n$ and $d$,
where the number of variables in the ring is $O(n)$
and the degree of the given generators of the ideal is $O(d)$.
Both $n$ and $d$ are positive integers.

Here is the definition of the Mayr-Meyer ideals:
let $n, d \ge 2$ be integers and $k$ a field.
Let
$s, f, s_{r+1}, f_{r+1},
b_{r1}, b_{r2}, b_{r3}, b_{r4}, c_{r1}, c_{r2}, c_{r3}, c_{r4}$
be variables over $k$,
with $r = 0, 1, \ldots, n-1$.
The notation here closely follows that of [K].
Set
$$
S = k[s = s_0, f = f_0, s_{r+1}, f_{r+1}, b_{ri}, c_{ri}| r = 0, \ldots, n-1; i = 1, \ldots, 4].
$$
Thus $S$ is a polynomial ring of dimension $10n+2$.
The following generators define the Mayr-Meyer ideal $J_l(n,d)$
(subscript $l$ for ``long",
there will be a ``shortened" version later on):
first the four level 0 generators:
$$
H_{0i} = c_{0i} \left(s -fb_{0i}^d\right), i = 1, 2, 3, 4;
$$
then the first six level $r$ generators, $r = 1, \ldots, n$:
$$
\eqalignno{
& H_{r1} = s_r - s_{r-1} c_{r-1,1}, \cr
& H_{r2} = f_r - s_{r-1} c_{r-1,4}, \cr
& H_{r3} = f_{r-1} c_{r-1,1} - s_{r-1} c_{r-1,2}, \cr
& H_{r4} = f_{r-1} c_{r-1,4} - s_{r-1} c_{r-1,3}, \cr
& H_{r5} = s_{r-1} \left( c_{r-1,3}-c_{r-1,2} \right), \cr
& H_{r6} = f_{r-1} \left( c_{r-1,2} b_{r-1,1}-c_{r-1,3} b_{r-1,4} \right), \cr
}
$$
the last four level $r$ generators, $r= 2, ..., n-1$:
$$
H_{r,6+i} = f_{r-1} c_{r-1,2} c_{ri}
\left( b_{r-1,2}-b_{ri} b_{r-1,3} \right),
i = 1, \ldots, 4,
$$
and the last level $n$ generator:
$$
H_{n7} = f_{n-1} c_{n-1,2}
\left( b_{n-1,2}-b_{n-1,3} \right).
$$
The maximum degree of a given generator of $J_l(n,d)$ is $d + 2$.
The degree 1 element $s_n - f_n$ of $S$ is in $J_l(n,d)$,
and when written as an $S$-linear combination of the given generators,
the $S$-coefficient of $H_{04}$ has degree
which is doubly exponential in $n$
(see [MM], [BS], [K]).

%The goal of this paper is to examine the associated primes of $J_l(n,d)$,
%and see if the number or the structure of these
%is what causes the doubly exponential behavior.  HERE too ambitious

The following summarizes the elementary facts used in the paper:

\vskip2ex\sectno=0\noindent{\bf Facts:}
\bgroup\parindent=3em
\fact
\label{\factsub}
For any ideals $I, I'$ and $I''$ with $I \subseteq I''$,
$(I + I') \cap I'' = I + I' \cap I''$.

\fact
For any ideal $I$ and element $x$,
$(x) \cap I = x (I : x)$.

\fact
\label{\factfewvars}
Let $x_1, \ldots, x_n$ be variables over a ring $R$.
Let $S = R[x_1, \ldots, x_n]$.
For any $f_1 \in R$,
$f_2 \in R[x_1]$,
$\ldots$, $f_n \in R[x_1, \ldots, x_{n-1}]$,
let $L$ be the ideal $(x_1 - f_1, \ldots, x_n - f_n)S$ in $S$.
Then an ideal $I$ in $R$ is primary (respectively, prime)
if and only if $IS + L$ is a primary (respectively, prime) in $S$.
Furthermore,
$\cap_i q_i = I$ is a primary decomposition of $I$
if and only if $\cap_i (q_iS + L)$ is a primary decomposition of $IS + L$.

\fact
\label{\factdec}
Let $x$ be an element of a ring $R$ and $I$ an ideal.
Suppose that there is an integer $k$ such that for all $m$,
$I : x^m \subseteq I : x^k$.
Then $I = \left(I : x^k\right) \cap \left(I + (x^k)\right)$.
Thus to find a (possibly redundant) primary decomposition of $I$
it suffices to find primary decompositions of (possibly larger) $I : x^k$
and of $I + (x^k)$.

\egroup

We immediately apply the last fact:
in order to find a primary decomposition of the Mayr-Meyer ideals $J_l(n,d)$,
by the structure of the $H_{r1}, H_{r2}$ and by Fact~\factfewvars,
it suffices to find a primary decomposition of the ideals $J(n,d)$
obtained from $J_l(n,d)$
by rewriting the variables $s_r, f_r$ in terms of other variables,
and then omitting the generators $H_{r1}, H_{r2}$, $r \ge 1$.
An ideal $q$ is a component (resp.\ associated prime) of $J(n,d)$
if and only if $(q + (H_{r1}, H_{r2} | r))S$ is
a component (resp.\ associated prime) of $J_l(n,d)$.
Thus to simplify the notation,
throughout we will be searching for
the primary components and associated primes
of the ``shortened'' Mayr-Meyer ideals $J(n,d)$
in a smaller polynomial ring $R$ obtained as above.
When we list the new generators explicitly,
the case $n = 1$ is rather special.
In fact,
the primary decomposition in the case $n = 1$ is very different
from the case $n \ge 2$,
and is given in [S1].
In this paper it is always assumed that $n \ge 2$.

Thus explicitly,
we will be working with the following ``shortened" Mayr-Meyer ideals:
for any fixed integers $n, d \ge 2$,
$R = k[s, f, b_{ri}, c_{ri}| r = 0, \ldots, n-1; i = 1, \ldots, 4]$,
a polynomial ring in $8n+2$ variables,
and $J(n,d)$ is the ideal in $R$
generated by the following polynomials $h_{ri}$:
first the four level 0 generators:
$$
h_{0i} = c_{0i} \left(s -fb_{0i}^d\right), i = 1, 2, 3, 4;
$$
then the eight level 1 generators:
$$
\eqalignno{
& h_{13} = fc_{01} - s c_{02}, \cr
& h_{14} = fc_{04} - s c_{03}, \cr
& h_{15} = s \left(c_{03} - c_{02} \right), \cr
& h_{16} = f \left(c_{02} b_{01} - c_{03} b_{04} \right), \cr
& h_{1,6+i} = f c_{02} c_{1i} \left(b_{02}-b_{1i} b_{03}\right), i = 1, \ldots, 4, \cr
}
$$
the first four level $r$ generators, $r = 2, \ldots, n$:
$$
\eqalignno{
&h_{r3} = s c_{01} c_{11} \cdots c_{r-3,1} \left(
c_{r-2,4} c_{r-1,1} - c_{r-2,1} c_{r-1,2}
\right), \cr
& h_{r4} = s c_{01} c_{11} \cdots c_{r-3,1} \left(
c_{r-2,4} c_{r-1,4} - c_{r-2,1} c_{r-1,3}
\right), \cr
& h_{r5} = s c_{01} c_{11} \cdots c_{r-2,1}
\left( c_{r-1,3}-c_{r-1,2} \right), \cr
& h_{r6} = s c_{01} c_{11} \cdots c_{r-3,1} c_{r-2,4}
\left( c_{r-1,2} b_{r-1,1}-c_{r-1,3} b_{r-1,4} \right), \cr
}
$$
the last four level $r$ generators, $r= 2, ..., n-1$:
$$
h_{r,6+i} = s c_{01} c_{11} \cdots c_{r-3,1} c_{r-2,4} c_{r-1,2} c_{ri}
\left( b_{r-1,2}-b_{ri} b_{r-1,3} \right),
i = 1, \ldots, 4,
$$
and the last level $n$ generator:
$$
h_{n7} = s c_{01} c_{11} \cdots c_{n-3,1} c_{n-2,4} c_{n-1,2}
\left( b_{n-1,2}-b_{n-1,3} \right).
$$
For simpler notation,
$J(n,d)$ will often be abbreviated to $J$.

Observe that the maximum degree of the given generators of $J(n,d)$
is $\max\{n+2, d+2\}$.
The image $s c_{01} c_{11} \cdots c_{n-2,1} (c_{n-1,1} - c_{n-1,4})$ of $s_n - f_n$
by construction lies in $J(n,d)$
and has degree $n+1$.
When this element is written as an $R$-linear combination of the $h_{ri}$,
the coefficient of $h_{04}$ is doubly exponential in $n$.
Note that the contrast between a number doubly exponential in $n$
and the degree $n+1$ of the input polynomial
arising from this instance of the ideal membership problem for $J(n,d)$
is not as striking as the contrast between a number doubly exponential in $n$
and the degree $1$ of the input polynomial
arising from the ideal membership example $s_n - f_n$ for $J_l(n,d)$.

Thus while $J(n,d)$ is a useful simplification of $J_l(n,d)$
as far as the primary decomposition and associated primes are concerned,
its doubly exponential nature is partially concealed.

This paper consists of two sections.
Section~\sectminimal\ is about all the minimal primes,
their components, and their heights.
For simplicity we assume that the underlying field $k$
is algebraically closed.
Then the number of minimal primes over $J(n,d)$ is $n(d')^2 + 20$
(Proposition~\propminprime),
where $d'$ is the largest divisor of $d$ which is relatively prime to
the characteristic of  the field.
All except 18 of the minimal components are simply the primes
(Proposition~\propmincomps).
Section \sectdblexpemb\ shows that
the doubly exponential behavior of the Mayr-Meyer ideals
is due to the existence of embedded primes.

The computation of embedded primes is tackled in [S2].
[S2] also constructs a new family of ideals
with the doubly exponential ideal membership problem.
Recursion can be applied to this new family
in the construction of the embedded prime ideals, see [S3].

\vskip 4ex

{\bf Acknowledgement.}
I thank Craig Huneke for suggesting this problem
all for all the conversations and enthusiasm for this research.

%\vfill\eject
\section{Minimal primes and their components}
\sectlabel{\sectminimal}
\pagelabel{\qsectminimal}

The minimal primes over $J(n,d)$
and their components are quite easy to compute.
Let $d'$ denote the largest divisor of $d$
which is relatively prime to the characteristic of the field.
Then there are are $n(d')^2 + 20$ minimal primes,
all but $18$ of which are their own primary components of $J(n,d)$.

The minimal primes are analyzed in two groups:
those on which $s$ and $f$ are non-zerodivisors,
and the rest of them.
The first group consists of $n(d')^2+1$ prime ideals.

The minimal primes not containing $sf$ are denoted $P_{r\underline{\hbox{\ \ }}}$,
where $r$ varies from $0$ to $n$,
and the other part $\underline{\hbox{\ \ }}$ of the subscript depends on $r$.
For the rest of the minimal primes the front part of the subscript
varies from $-1$ to $-4$.

\lemma
\label{\lmnosf}
Let $P$ be an ideal of $R$ containing $J$ such that
$s$ and $f$ are non-zerodivisors modulo $P$
(in particular $sf \not \in P$).
Let $r \in \{0, \ldots, n-1\}$.
Suppose that for all $j < r$ and all $i = 1, 2, 3, 4$,
$c_{ji}$ is not a zero-divisor modulo $P$.
Then
\itemitem{(1)}
For all $j \in \{0, \ldots, r\}$,
$$
c_{j3} - c_{j2}, c_{j4} - c_{j1}, c_{01}-c_{02}b_{01}^d \in P,
$$
and if $j > 0$,
$$
c_{j2} - c_{j1} \in P.
$$
\itemitem{(2)}
If $r > 0$,
$c_{ri} \in P$ for some $i \in \{1, 2, 3, 4\}$
if and only if
$c_{ri} \in P$ for all $i \in \{1, 2, 3, 4\}$.
\itemitem{(3)}
For all $j \in \{0, \ldots, r-1\}$,
$$
b_{j4} - b_{j1} \in P.
$$
Also, for all $j \in \{0, \ldots, r-2\}$,
$$
b_{j2}-b_{j+1,i}b_{j3} \in P, i = 1, 2, 3, 4.
$$
\itemitem{(4)}
Assume that $r > 0$.
Then for all $i, j \in \{1, 2, 3, 4\}$,
$$
b_{0i}^d - b_{0j}^d \in P.
$$
\itemitem{(5)}
Assume that $r > 0$
and that $P$ is a primary ideal such that no $b_{0i}$ lies in $\sqrt P$.
Then $s - fb_{01}^d \in P$,
and whenever $1 \le i < j \le 4$,
there exists a $(d')$th root of unity $\alpha_{ij} \in k$
such that $b_{0i} - \alpha_{ij} b_{0j} \in P$ and
$$
\alpha_{14} = 1,
\alpha_{24} = \alpha_{12}^{-1},
\alpha_{34} = \alpha_{13}^{-1}.
$$
\endb

\proof
By the assumption that $sf$ is a non-zerodivisor modulo $J$,
if $j = 0$,
$h_{15} = s \left(c_{03} - c_{02} \right)$ being in $P$
implies that $c_{03} - c_{02}$ is in $P$.
Also,
$h_{14}-h_{13}$ equals $f \left(c_{04} - c_{01} \right)$,
so that $c_{04} - c_{01} \in P$.
Note that $h_{01}-b_{01}^d h_{13} = s(c_{01}-c_{02}b_{01}^d)$,
so that $c_{01}-c_{02}b_{01}^d \in P$.
This proves (1) for $j = 0$.

Now assume that $j > 0$.
If $j \le r < n$,
$h_{j+1,5} = s c_{01} c_{11} \cdots c_{j-1,1}
\left( c_{j3}-c_{j2} \right)$ being in $P$
implies that $c_{j3}-c_{j2}$ is in $P$.
Furthermore,
$h_{j+1,4}-h_{j+1,3} = s c_{01} c_{11} \cdots c_{j-2,1} c_{j-1,4}
\left( c_{j4}-c_{j1} \right)$
so that $c_{j4}-c_{j1}$ is in $P$.
Then $h_{j+1,3}$ equals
$sc_{01} c_{11} \cdots c_{j-1,1} \left(c_{j1} - c_{j2}\right)$
modulo $(c_{j-1,4}-c_{j-1,1})$,
so that $c_{j1} - c_{j2}$ lies in $P$.
This proves (1).

With (1) established, (2) is an easy consequence.

To prove (3),
observe that modulo $(c_{03}-c_{02}) \subseteq P$,
$h_{16}$ equals $fc_{02} \left( b_{01} - b_{04}\right)$.
Hence if $r > 0$,
$b_{01} - b_{04}$ is in $P$.
If $0 \le j < r$,
$$
h_{j+1,6} \equiv s c_{01} c_{11} \cdots c_{j-2,1} c_{j-1,4} c_{j2}
\left(b_{j1} - b_{j4}\right)
\mod (c_{j3}-c_{j2}),
$$
hence $b_{j1} - b_{j4}$ is in $P$.
Furthermore,
for all $i = 1, \ldots, 4$,
$$
\eqalignno{
h_{1,6+i} &= f c_{02} c_{14}
\left( b_{02}-b_{1i} b_{03} \right) \in P, \cr
h_{j,6+i} &= s c_{01} c_{11} \cdots c_{j-3,1} c_{j-2,4} c_{j-1,2} c_{ji}
\left( b_{j-1,2}-b_{ji} b_{j-1,3} \right) \in P
\hbox{\ \ for $j > 1$},
}
$$
so that $b_{j-1,2}-b_{ji} b_{j-1,3}$ is in $P$
for all $j = 1, \ldots, r-1$ and all $i = 1, \ldots, 4$.
This proves (3).

If $r > 0$,
$h_{0i} = c_{0i} (s -fb_{0i}^d) \in P$ implies that $s -fb_{0i}^d \in P$.
Hence whenever $1 \le i < j \le 4$,
$f \left(b_{0i}^d-b_{0j}^d \right)$ is in $P$
so that $b_{0i}^d-b_{0j}^d$ is in $P$.
This proves (4),
and then (5) follows easily.
\qed

For notational purposes define the following ideals in $R$:
$$
\eqalignno{
\E &= (s - fb_{01}^d)
+ \left(b_{01} - b_{04}, b_{02}^d - b_{03}^d, b_{01}^d - b_{02}^d\right), \cr
\F &= \left( b_{02} -b_{11}b_{03}, b_{14}-b_{11},b_{13}-b_{11},
b_{12}-b_{11},b_{12}^d - 1\right) \cr
\C_r &= (c_{r1}, c_{r2}, c_{r3}, c_{r4}), r = 0, \ldots, n-1 \cr
\C_n &= (0), \cr
\D_0 &=
\left(c_{04}-c_{01},c_{03}-c_{02},c_{01}-c_{02}b_{01}^d \right), \cr
\D_r &= \left(c_{r4}-c_{r1},c_{r3}-c_{r2},c_{r2}-c_{r1} \right), r = 1, \ldots, n-1, \cr
\D_n &= (0), \cr
B_{0} &= B_{1} = (0), \cr
B_{r} &= \left(1-b_{2i}, 1-b_{3i}, \ldots, 1-b_{ri} |
i = 1, \ldots, 4 \right), r = 2, \ldots, n-1. \cr
%B_{kr} &= \left(1-b_{ki}, 1-b_{k+1,i}, \ldots, 1-b_{ri} |
%i = 1, \ldots, 4 \right), r = 2, \ldots, n-1. \cr
%\D_n &= (0). \cr
}
$$

With the previous lemma and this notation then:

\prop
\label{\propnosf}
Let $P$ be a minimal prime ideal containing $J$ and not containing $sf$.
\itemitem{(1)}
If $P$ contains one of the $c_{0i}$,
then $P$ equals the height four prime ideal
$$
P_0 = (c_{01},c_{02},c_{03},c_{04}) = \C_0.
$$
\itemitem{(2)}
If $P$ contains no $c_{ji}$,
set $r = n$,
otherwise
set $r$ to be the smallest integer such that $P$ contains some $c_{ri}$.
If $r = 1$,
$P$ contains
$$
p_1 = \C_1 + \E + \D_0,
$$
and if $r > 1$,
$P$ contains
$$
p_r = \C_r + \E + \F + B_{r-1} + \D_0 + \D_1 + \cdots + \D_{r-1}.
$$
\itemitem{(3)}
For all $r = 1, \ldots, n$,
$J \subseteq p_r$.
\endb

\proof
Suppose that $P$ contains $c_{01}$ or $c_{04}$.
Then by Lemma~\lmnosf, $P$ contains $c_{02}b_{0i}^d$.
If $c_{02}$ is not in $P$,
then $b_{02} \in P$,
hence as $h_{02} = c_{02} (s -f b_{02}^d) \in P$,
necessarily $c_{02} s \in P$,
contradicting the choice.
Thus necessarily $P$ contains $c_{02}$.
Then by Lemma~\lmnosf, $P$ contains all the $c_{0i}$.
As $P_0$ contains $J$,
this verifies (1).

If $r \ge 1$,
$p_r$ obviously contains $J$, thus verifying (3).
If $b_{02}$ is in $P$,
then as above also $c_{02} s$ is in $P$,
contradicting the assumptions.
Thus $b_{02}$ is not in $P$
and (2) follows for the case $r = 1$ by Lemma~\lmnosf.
Now let $r > 1$.
By Lemma~\lmnosf,
it remains to prove that $\F + B_{r-1} \subseteq P$
when $r > 1$.
As $b_{j2}-b_{j+1,i}b_{j3} \in P$
for all $j = 0, \ldots, r-2$, $i = 1, \ldots, 4$,
it follows that $(b_{j+1,i}-b_{j+1,i'})b_{j3}$ is in $P$
for any $i, i' \in \{1, 2, 3, 4\}$.
If $b_{j3} \in P$,
by an application of Lemma~\lmnosf~(3),
$b_{j-1,2} \in P$,
whence $b_{j-2,2} \in P$,
$\ldots$,
$b_{02}$ is in $P$,
which is a contradiction.
Thus necessarily
$b_{j+1,i}-b_{j+1,i'}$ is in $P$ for all $j = 0, \ldots, r-2$,
or that
$b_{j-1,i}-b_{j-1,i'}$ is in $P$ for all $j = 2, \ldots, r$.
Once this is established,
then $h_{j,6+i}$ equals
$s c_{01} c_{11} \cdots c_{j-3,1} c_{j-2,4} c_{j-1,2} c_{ji} b_{j-1,3}
\left( 1-b_{ji} \right)$
modulo $P$
so that $1-b_{ji}$ is in $P$ for all $i = 1, \ldots, 4$
and all $j = 2, \ldots, r-1$.
A similar argument shows that $b_{11}^d - 1$ is in $P$.

The remaining case $r = n$ has essentially the same proof.
\qed

From this one can read off the minimal primes and components:

\prop
\label{\propprimesnosf}
All the minimal prime ideals over $J$
which do not contain $sf$ are
$$
\eqalignno{
&P_0, \cr
&P_{1\alpha\beta} = p_1 + (b_{01}-\alpha b_{02},b_{02}-\beta b_{03}), \cr
&P_{r\alpha\beta} = p_r + (b_{01}-\alpha b_{02},b_{02}-\beta b_{03},
\beta -b_{1i}| i = 1, \ldots, 4),\cr
}
$$
where $\alpha$ and $\beta$ vary over the $(d')$th roots of unity.
The heights of these ideals are as follows:
$\hht(P_0) = 4$,
for $r \in \{1, \ldots, n-1\}$,
$\hht(P_{r\alpha\beta}) = 7r + 4$,
and $\hht(P_{n\alpha\beta}) = 7n$.
The components of $J(n,d)$ corresponding to these prime ideals
are the primes themselves.

Furthermore,
with notation as in the previous proposition,
for all $r \ge 1$,
$\cap_{\alpha,\beta} P_{r\alpha\beta} = p_r$.
\endb

\proof
The case of $P_0$ is trivial.
It is easy to see that for $r > 0$,
the listed primes $P_{r\alpha\beta}$ are minimal over $p_r$
and that the intersection of the $(d')^2$ $P_{r\alpha\beta}$ equals $p_r$.
It is trivial to calculate the heights,
and it is straightforward to prove the last statement.
\qed

This completes the list of all the minimal primes over $J(n,d)$
which do not contain $s$ and $f$.
The next group of minimal primes all contain $s$:

\prop
\label{\propprimes}
Let $P$ be a prime ideal minimal over $J$.
If $P$ contains $s$,
then $P$ is one of the following 19 prime ideals:
$$
\eqalignno{
P_{-1} &= (s,f), \cr
P_{-2} &= (s, c_{01}, c_{02}, c_{04}, b_{03}, b_{04}) \cr
P_{-3} &= (s, c_{01}, c_{04}, b_{02}, b_{03},c_{02} b_{01} - c_{03} b_{04}), \cr
P_{-4\Lambda} &= \left(c_{1i} | i \not \in \Lambda \right)
+ \left(b_{1i} | i \in \Lambda \right)
+ \left(s, c_{01}, c_{03}, c_{04}, b_{01}, b_{02} \right), \cr
}
$$
as $\Lambda$ varies over the subsets of $\{1, 2, 3, 4\}$.
The heights of these primes are $2, 6, 6$ and $10$, respectively.
\endb

\proof
Note that
$$
J + (s) =
\left(c_{0i} fb_{0i}^d, fc_{02}c_{1i}\left(b_{02}-b_{1i}b_{03}\right)
| i = 1, 2, 3, 4\right)
+ \left(s, fc_{01}, fc_{04}, f \left(c_{02} b_{01} - c_{03} b_{04} \right)\right).
$$
If $P$ contains $f$, it certainly equals $P_{-1}$.
Now assume that $P$ does not contain $f$.
Then $P$ is minimal over
$$
\left(c_{0i} b_{0i}^d, c_{02}c_{1i}\left(b_{02}-b_{1i}b_{03}\right)
| i = 1, 2, 3, 4\right)
+ \left(s, c_{01}, c_{04}, c_{02} b_{01} - c_{03} b_{04} \right).
$$
If $c_{02} \in P$,
then $P$ is minimal over
$$
\left(c_{03} b_{03}^d, s, c_{01}, c_{02}, c_{04}, c_{03} b_{04} \right),
$$
so it is either
$\left(s, c_{01}, c_{02}, c_{03}, c_{04}\right)$
or
$\left(s, c_{01}, c_{02}, c_{04}, b_{03}, b_{04} \right)= P_{-2}$.
However,
the first option is not minimal over $J$
as it strictly contains $P_0$ from Proposition~\propnosf.

Now assume that $P$ does not contain $fc_{02}$.
Then $P$ is minimal over
$$
(b_{02}, c_{03} b_{03}^d)
+ \left(c_{1i}b_{1i}b_{03} | i = 1, 2, 3, 4\right)
+ \left(s, c_{01}, c_{04}, c_{02} b_{01} - c_{03} b_{04} \right).
$$
If $P$ contains $b_{03}$,
then $P = (s, c_{01}, c_{04}, b_{02}, b_{03},c_{02} b_{01} - c_{03} b_{04})$,
which is $P_{-3}$.

Finally,
assume that $P$ does not contain $fc_{02}b_{03}$.
Then $P$ is minimal over
$$
(b_{02}, c_{03})
+ \left(c_{1i}b_{1i}| i = 1, 2, 3, 4\right)
+ \left(s, c_{01}, c_{04}, b_{01} \right),
$$
whence $P$ is one of the $P_{-4\Lambda}$.
\qed

It turns out that there are no other minimal primes over $J(n,d)$:

\prop
\label{\propminprime}
The prime ideals from the previous two propositions
are the only prime ideals minimal over $J$.
Thus there are $1 + n(d')^2 + 3 + 2^4 = n(d')^2 + 20$ minimal primes.
\endb

\proof
Proposition~\propprimesnosf\ determined all the minimal primes over $J$
not containing $sf$,
and Proposition~\propprimes\ determined all those minimal primes
which contain $s$.
It remains to find all the prime ideals containing $f$ and $J$ but not $s$.
As
$J + (f)$ contains
$\left(c_{0i} s | i = 1, 2, 3, 4\right)$,
a minimal prime ideal containing $J + (f)$ but not $s$
contains, and even equals $(f, c_{01},c_{02},c_{03},c_{04})$.
However,
this prime ideal properly contains $P_0$,
and hence is not minimal over $J$.
The proposition follows as there are no containment relations
among the given primes.
\qed

The $n(d')^2 + 20$ minimal primary components can be easily computed:

\prop
\label{\propmincomps}
For all possible subscripts $\circ$,
let $p_{\circ}$ be the $P_{\circ}$-primary component of $J$.
Then
$$
\eqalignno{
p_{-2} &= (s, c_{01}, c_{02}, c_{04}, b_{03}^d, b_{04}), \cr
%p_{-3} &= (s, c_{01}, c_{04}, b_{02}, b_{03}, c_{02}b_{01} - c_{03} b_{04}), \cr
p_{-4\Lambda} &= \left(c_{1i} | i \not \in \Lambda \right)
+ \left(b_{1i}^d, b_{02}-b_{1i}b_{03},b_{1i}-b_{1j} | i, j \in \Lambda \right)
+ \left(s, c_{01}, c_{03}, c_{04}, b_{01}, b_{02}^d \right), \cr
}
$$
and all the other $p_{\circ}$ equal $P_{\circ}$.
\endb

\proof
By Proposition~\propprimesnosf,
it remains to calculate $p_{-1},p_{-2},p_{-3}$ and $p_{-4\Lambda}$.
As $c_{03}-c_{02}$
is not an element of $P_{-1}, P_{-2}, P_{-3}$ and $P_{-4\Lambda}$,
and since $h_{15} = s \left(c_{03} - c_{02} \right)$ is in $J$,
it follows that $s \in p_{-1}, p_{-2}, p_{-3}$ and $p_{-4\Lambda}$.
Then $c_{01} f b_{01}^d \in p_{-1}$,
so that $f \in p_{-1}$,
and so $p_{-1} = P_{-1}$.

As $h_{13} = fc_{01} - s c_{02},
h_{14} = fc_{04} - s c_{03}$ are in $J$,
then
$fc_{01}, fc_{04} \in p_{-2}, p_{-3}$ and $p_{-4\Lambda}$,
whence
$c_{01}, c_{04} \in p_{-2}, p_{-3}, p_{-4\Lambda}$.
For all $i = 1, \ldots, 4$,
as $h_{1,6+i} = fc_{02}c_{1i}(b_{02}-b_{1i}b_{03}) \in J$,
it follows that $b_{02}-b_{1i}b_{03} \in p_{-3}$.
Thus as $b_{11}-b_{12} \not \in P_{-3}$,
it follows that $b_{03}$ and hence also $b_{02}$ are in $p_{-3}$.
Now it is clear that $p_{-3}$ is the $P_{-3}$-primary component of $J$.

Further for $i = 2, 3$,
$c_{0i} f b_{0i}^d \in p_j$ implies that
$b_{02}^d \in p_{-4\Lambda}$,
$b_{03}^d \in p_{-2}$,
$c_{02} \in p_{-2}$,
and
$c_{03} \in p_{-4\Lambda}$.
As $h_{16} = f \left(c_{02} b_{01} - c_{03} b_{04} \right)$ is in $J$,
then
$f c_{03} b_{04}$ is in $p_{-2}$ so that $b_{04}$ is in $p_{-2}$.
Also
$f c_{02} b_{01}$ is in $p_{-4\Lambda}$ so that $b_{01}$ is in $p_{-4\Lambda}$.
Thus the $P_{-2}$-primary component contains $p_{-2}$.
But $p_{-2}$ contains $J$,
so $p_{-2}$ is the $P_{-2}$-primary component of $J$.

Lastly,
as $J$ contains $h_{1,6+i}$, $i = 1, \ldots, 4$,
each $p_{-4\Lambda}$ contains each $c_{1i}(b_{02}-b_{1i}b_{03})$.
If $i \not \in \Lambda$,
then $b_{02}-b_{1i}b_{03}$ is not in $P_{-4\Lambda}$,
so that $c_{1i} \in p_{-4\Lambda}$.
If instead $i \in \Lambda$,
then $c_{1i} \not \in P_{-4\Lambda}$,
so that
$b_{02}-b_{1i}b_{03}$ is in $p_{-4\Lambda}$.
Hence $b_{02}^d-b_{1i}^db_{03}^d$ is in $p_{-4\Lambda}$,
so that as $b_{02}^d \in p_{-4\Lambda}$,
so is $b_{1i}^db_{03}^d$.
Hence
$b_{1i}^d$ is in $p_{-4\Lambda}$.
Furthermore,
for $i, j \in \Lambda$,
$b_{03}(b_{1j}-b_{1i}) = (b_{02}-b_{1i}b_{03})- (b_{02}-b_{1j}b_{03})$
is in $p_{-4\Lambda}$,
so that
$b_{1j}-b_{1i}$ is in $p_{-4\Lambda}$.
Thus
$$
p_{-4\Lambda} \supseteq
\left(c_{1i} | i \not \in \Lambda \right)
+ \left(b_{1i}^d, b_{02}-b_{1i}b_{03},b_{1i}-b_{1j} | i, j \in \Lambda \right)
+ \left(s, c_{01}, c_{03}, c_{04}, b_{01}, b_{02}^d \right),
$$
but the latter ideal is primary and contains $J$,
so equality holds.
\qed

The structure of $p_{-2}$ says that

\prop
\label{\propnotrad}
For $n, d \ge 2$,
$J(n,d)$ is not a radical ideal.
\qed
\endb

Here is the table of all the minimal primes over $J(n,d)$,
where $\alpha$ and $\beta$ are $(d')$th roots of unity:

\vskip 2ex
\halign{\mvrule \hskip 0.1em #\hfil \hskip 0.1em
&& \mvrule \hskip 0.5em \relax #\hfil \mvrule \cr
\noalign{\hrule}
%\noalign{\vskip 1pt \hrule \vskip 1pt}
minimal prime & height & component of $J(n,d)$ \cr
\noalign{\hrule}
\noalign{\hrule}
$P_0 = (c_{01}, c_{02}, c_{03}, c_{04})$ & $4$ & $p_0 = P_0$ \cr
\noalign{\hrule}
$P_{1\alpha\beta}= p_1 + (b_{01}- \alpha b_{02}, b_{02}-\beta b_{03})$
& $11$ & $p_{1\alpha\beta}= P_{1\alpha\beta}$ \cr
\noalign{\hrule}
$P_{r\alpha\beta}=p_r$ & $7r + 4$
& $p_{r\alpha\beta}= P_{r\alpha\beta}$, $2 \le r < n$ \cr
\hskip 1em $+ (b_{01}- \alpha b_{02}, b_{02}-\beta b_{03}, \beta-b_{1i})$
& $7n$ & $p_{r\alpha\beta}= P_{r\alpha\beta}$, $r =n$ \cr
\noalign{\hrule}
$P_{-1}= (s,f)$ & $2$ & $p_{-1} = P_{-1}$ \cr
\noalign{\hrule}
$P_{-2}= (s, c_{01}, c_{02}, c_{04}, b_{03}, b_{04})$ & $6$
& \ $p_{-2}= (s, c_{01}, c_{02}, c_{04}, b_{03}^d, b_{04})$ \cr
\noalign{\hrule}
$P_{-3} = (s, c_{01}, c_{04}, b_{02}, b_{03}, c_{02}b_{01}-c_{03}b_{04})$
& $6$ & $p_{-3} = P_{-3}$ \cr
\noalign{\hrule}
$P_{-4\Lambda} =
(s, c_{01}, c_{03}, c_{04}, b_{01}, b_{02})$ & $10$
& $p_{-4\Lambda} =(s, c_{01}, c_{03}, c_{04}, b_{01}, b_{02}^d)$\cr
$\hskip1em + (c_{1i}, b_{1j} |i \not \in \Lambda, j \in \Lambda)$ &
& $+ (c_{1i}|i \not \in \Lambda)$ \cr
& & $+ (b_{1j}^d, b_{02}-b_{1j}b_{03}, b_{1j}-b_{1j'}
|j,j' \in \Lambda)$ \cr
\noalign{\hrule}
}

\vskip 2ex

%\vfill\eject
\section{Doubly exponential behavior is due to embedded primes}
\sectlabel{\sectdblexpemb}
\pagelabel{\qsectdblexpemb}

In this section we compute the intersection of all the minimal components of $J(n,d)$,
and show that the element $sc_{01} \cdots c_{n-2,1} (c_{n-1,1} - c_{n-1,4})$,
which gave the doubly exponential membership property for $J(n,d)$,
does not give the doubly exponential membership property for the intersection
of the minimal components.
This proves that the doubly exponential behavior of the Mayr-Meyer ideals
is due to the existence of embedded primes.

%The computation of the intersection of the minimal components is rather lengthy and laborious.

First define the ideal
$$
p_{-4}= (s, c_{01}, c_{03}, c_{04}, b_{01}, b_{02}^d)
+ (c_{1i}(b_{02}-b_{1i}b_{03}), c_{1i} b_{1i}^d,
c_{1i}c_{1j}(b_{1i}-b_{1j})| i,j = 1, \ldots, 4).
$$
Note that $p_{-4}:c_{11} = p_{-4}:c_{11}^2$,
so that by Fact~\factdec,
$p_{-4} = (p_{-4}:c_{11}) \cap (p_{-4} + (c_{11}))$.
\pagelabel{\qcalc}
Similarly,
$$
\eqalignno{
p_{-4} &= (p_{-4}:c_{11} c_{12}) \cap((p_{-4}:c_{11})+( c_{12}))
\cap ((p_{-4} + (c_{11})) : c_{12})
\cap (p_{-4} + (c_{11} ,c_{12})) \cr
&= \cdots \cr
&= \bigcap_{\Lambda}
\left( \left( \left( \left(p_{-4} *_1 c_{11}\right) *_2 c_{12}\right)
*_3 c_{13}\right)*_4 c_{14}\right), \cr
}
$$
where $*_i$ vary over the operations colon and addition.
But the resulting component ideals are just the various $p_{-4\Lambda}$,
so that
$$
p_{-4} = \bigcap_{\Lambda} p_{-4\Lambda}.
%\left( \left(c_{1i} | i \not \in \Lambda \right)
%+ \left(b_{1i}^d, b_{02}-b_{1i}b_{03},b_{1i}-b_{1j} | i, j \in \Lambda \right)
%+ \left(s, c_{01}, c_{03}, c_{04}, b_{01}, b_{02}^d \right) \right).
$$
%is simply the intersection of all the $p_{-4\Lambda}$.
Next we compute the intersection of $p_{-4}$ and $p_{-2}$
(using Fact~\factsub):
$$
\eqalignno{
p_{-2} &\cap p_{-4}
=
(s, c_{01}, c_{04}) + (c_{02}, b_{03}^d, b_{04}) \cap p_{-4} \cr
&=
(s, c_{01}, c_{04}) + (c_{02}, b_{03}^d, b_{04}) \cdot
(c_{03}, b_{01}, b_{02}^d,c_{1i}(b_{02}-b_{1i}b_{03}), c_{1i} b_{1i}^d,
c_{1i}c_{1j}(b_{1i}-b_{1j})), \cr
}
$$
so that
$$
\eqalignno{
p_{-2} &\cap p_{-3} \cap p_{-4}
=
(s, c_{01}, c_{04}, c_{02} b_{01} - c_{03} b_{04}) \cr
&\hskip2em
+ (b_{03}^d) 
(c_{03}, b_{01}, b_{02}^d,c_{1i}(b_{02}-b_{1i}b_{03}), c_{1i} b_{1i}^d,
c_{1i}c_{1j}(b_{1i}-b_{1j})) \cr
&\hskip2em
+ (c_{02}, b_{04}) \cdot
(b_{02}^d,c_{1i}(b_{02}-b_{1i}b_{03})) \cr
&\hskip2em
+ (c_{02}, b_{04}) \cdot
(c_{03}, b_{01}, c_{1i} b_{1i}^d,c_{1i}c_{1j}(b_{1i}-b_{1j})) \cdot (b_{02}, b_{03}) \cr
&=
(s, c_{01}, c_{04}, c_{02} b_{01} - c_{03} b_{04}) \cr
&\hskip2em
+ (b_{03}^d) 
(c_{03}, b_{01}, b_{02}^d,c_{1i}(b_{02}-b_{1i}b_{03}), c_{1i} b_{1i}^d,
c_{1i}c_{1j}(b_{1i}-b_{1j})) \cr
&\hskip2em
+ (c_{02}, b_{04})
(b_{02}^d,c_{1i}(b_{02}-b_{1i}b_{03}),c_{03}b_{02}, c_{03}b_{03},
b_{01}b_{02},b_{01}b_{03},
c_{1i} b_{1i}^db_{02},c_{1i} b_{1i}^db_{03}). \cr
}
$$
Thus the intersection of the minimal components of $J(n,d)$
which contain $s$ equals:
$$
\eqalignno{
p_{-1} &\cap p_{-2} \cap p_{-3} \cap p_{-4}
=
(s, fc_{01}, fc_{04}, f(c_{02} b_{01} - c_{03} b_{04})) \cr
&\hskip2em
+ f(b_{03}^d) 
(c_{03}, b_{01}, b_{02}^d,c_{1i}(b_{02}-b_{1i}b_{03}), c_{1i} b_{1i}^d,
c_{1i}c_{1j}(b_{1i}-b_{1j})) \cr
&\hskip2em
+ f(c_{02}, b_{04})
(b_{02}^d,c_{1i}(b_{02}-b_{1i}b_{03}),c_{03}b_{02}, c_{03}b_{03},
b_{01}b_{02},b_{01}b_{03},
c_{1i} b_{1i}^db_{02},c_{1i} b_{1i}^db_{03}) \cr
&=
J + (s)
+ f b_{03}^d 
(b_{01}, b_{02}^d,c_{1i}(b_{02}-b_{1i}b_{03}), c_{1i} b_{1i}^d,
c_{1i}c_{1j}(b_{1i}-b_{1j})) \cr
&\hskip2em
+ f c_{02}
(c_{03}b_{02}, c_{03}b_{03},
b_{01}b_{02},b_{01}b_{03},c_{1i} b_{1i}^db_{03}) \cr
&\hskip2em
+ fb_{04}
(b_{02}^d,c_{1i}(b_{02}-b_{1i}b_{03}),
b_{01}b_{02},b_{01}b_{03},
c_{1i} b_{1i}^db_{02},c_{1i} b_{1i}^db_{03}). \cr
}
$$
We can simplify this intersection in terms of the generators of $J$
if we first intersect the intersection with the minimal component $p_0$:
$$
\eqalignno{
p_{0} &\cap \cdots \cap p_{-4}
=
J + s C_0
+ f b_{03}^d C_0 
(b_{01}, b_{02}^d,c_{1i}(b_{02}-b_{1i}b_{03}), c_{1i} b_{1i}^d,
c_{1i}c_{1j}(b_{1i}-b_{1j})) \cr
&\hskip2em
+ f c_{02}
(c_{03}b_{02}, c_{03}b_{03},
b_{01}b_{02},b_{01}b_{03},c_{1i} b_{1i}^db_{03}) \cr
&\hskip2em
+ fb_{04} \C_0
(b_{02}^d,c_{1i}(b_{02}-b_{1i}b_{03}),
b_{01}b_{02},b_{01}b_{03},
c_{1i} b_{1i}^db_{02},c_{1i} b_{1i}^db_{03}). \cr
}
$$
As $f \C_0 \subseteq J + f(c_{02}, c_{03})$
and $s \C_0 \subseteq J + s \D_0 + (fc_{02}b_{02}^d)$,
it follows that
$$
\eqalignno{
p_{0} \cap \cdots \cap p_{-4}
&=
J + s C_0
+ f b_{03}^d (c_{02}, c_{03}) 
(b_{01}, b_{02}^d,c_{1i}(b_{02}-b_{1i}b_{03}), c_{1i} b_{1i}^d,
c_{1i}c_{1j}(b_{1i}-b_{1j})) \cr
&\hskip2em
+ f c_{02}
(c_{03}b_{02}, c_{03}b_{03},
b_{01}b_{02},b_{01}b_{03},c_{1i} b_{1i}^db_{03}) \cr
&\hskip2em
+ fb_{04} (c_{02}, c_{03})
(b_{02}^d,c_{1i}(b_{02}-b_{1i}b_{03}),
b_{01}b_{02},b_{01}b_{03},
c_{1i} b_{1i}^db_{02},c_{1i} b_{1i}^db_{03}) \cr
&=J + s \C_0
+ f c_{02}(c_{03}b_{02}, c_{03}b_{03},
b_{01}b_{02},b_{01}b_{03},c_{1i} b_{1i}^db_{03}) \cr
&=J + s \D_0
+ f c_{02}(c_{03}b_{02}, c_{03}b_{03},
b_{01}b_{02},b_{01}b_{03},c_{1i} b_{1i}^db_{03},b_{02}^d). \cr
}
$$

Next we compute the intersection of all the minimal components of $J(n,d)$
which do not contain $s$:

\lemma
For $2 \le r \le n$,
$$
p_1 \cap p_2 \cap \cdots \cap p_r = E + \D_0 + \C_1 F
+ \sum_{i=0}^{r-1} \C_1 \C_2 \cdots \C_i \left(\D_{i+1} + B_i \right)
+ \C_1 \C_2 \cdots \C_r.
$$
\endb

\proof
When $r = 2$,
$$
\eqalignno{
p_1 \cap p_2 &=
\left(\C_1 + E + \D_0\right)
\cap \left(\C_2 + E + F + \D_0 + \D_1 + B_1\right) \cr
&=
E + \D_0 +
\C_1 \cap \left(\C_2 + E + F + \D_0 + \D_1 + B_1\right) \cr
&=
E + \D_0 + \D_1 +
\C_1 \cap \left(\C_2 + E + F + \D_0 + B_1\right) \cr
&=
E + \D_0 + \D_1 +
\C_1 \cdot \left(\C_2 + E + F + \D_0 + B_1\right) \cr
&=
E + \D_0 + \D_1 + \C_1 F +
\C_1 \cdot \left(\C_2 + B_1\right), \cr
}
$$
which starts the induction.
Then by induction assumption for $r \ge 2$ and $r \le n-1$,
$$
\eqalignno{
p_1 \cap \cdots \cap p_{r+1} &=
\left(E + \D_0 + \C_1 F
+ \sum_{i=0}^{r-1} \C_1 \C_2 \cdots \C_i \left(\D_{i+1} + B_i \right)
+ \C_1 \cdots \C_r \right)
\cap p_{r+1} \cr
&=
E + \D_0 + \C_1 F
+ \sum_{i=0}^{r-1} \C_1 \C_2 \cdots \C_i \left(\D_{i+1} + B_i \right)
+ (\C_1 \cdots \C_r) \cap p_{r+1}, \cr
}
$$
and by multihomogeneity,
the last intersection equals
$$
\C_1 \cdots \C_r \left(\C_{r+1} + E + F + B_r \right)
+
\sum_{j = 1}^r ( \D_j \prod_{k \not = j} ^r \C_k).
$$
Combining the last two displays proves the lemma.
\qed

%\remark
%One can prove that $J : s^2 \subsetneq J : s^3 = J : s^4$.
%Thus by Fact~\factdec,
%a primary decomposition of $J$ can be obtained from
%$J = (J : s^3) \cap (J + (s^3))$.
%The first component is easy:
%$J : s^3$ equals the intersection of all the $p_i$, $i \ge 0$.
%The proof takes 6 pages
%and is included in the ``omitted'' paper found on my web site
%{\tt http://math.nmsu.edu/\~{}iswanson}.
%Unfortunately,
%the decomposition of $J + (s^3)$ has eluded me.
%See the paper on my web page for the various attempts at finding
%the embedded primes of $J(n,d)$.
%\endb

Thus the intersection of all the minimal components of $J(n,d)$ equals:
$$
\eqalignno{
\bigcap_{r = -4} ^n p_r &=
\left(J + s \D_0 + f c_{02}
(c_{03}b_{02}, c_{03}b_{03},
b_{01}b_{02},b_{01}b_{03},c_{1i} b_{1i}^db_{03},b_{02}^d) \right)
\cap \bigcap_{r = 1} ^n p_r \cr
&= J + s \D_0
+ \left(f c_{02}
(c_{03}b_{02}, c_{03}b_{03},
b_{01}b_{02},b_{01}b_{03},c_{1i} b_{1i}^db_{03},b_{02}^d) \right) 
\cap \bigcap_{r = 1} ^n p_r \cr
&= J + s \D_0
+ f c_{02}\left(
(c_{03}b_{02}, c_{03}b_{03},
b_{01}b_{02},b_{01}b_{03},c_{1i} b_{1i}^db_{03},b_{02}^d)
\cap \bigcap_{r = 1} ^n p_r \right). \cr
}
$$
Let $A = (c_{03}b_{02}, c_{03}b_{03},
b_{01}b_{02},b_{01}b_{03},c_{1i} b_{1i}^db_{03},b_{02}^d)
\cap \bigcap_{r = 1} ^n p_r$.
Thus the intersection of all the minimal components of $J(n,d)$
equals $J + s\D_0 + fc_{02} A$.
Finding the generators of $A$ takes up most of the rest of this section.
We will use the decomposition 
$$
\eqalignno{
A &= (c_{03},b_{01},c_{1i} b_{1i}^d,b_{02}^d) \cap
(c_{03},b_{01},b_{03},b_{02}^d) \cap (b_{02}, b_{03})\cap \bigcap_{r = 1} ^n p_r, \cr
}
$$
and start computing $A$ via the indicated partial intersections,
again using Fact~\factsub:
$$
\eqalignno{
(b_{02}, b_{03}) \cap \bigcap_{r = 1} ^n p_r &=
(b_{02}^d - b_{03}^d, c_{11}(b_{02} - b_{11} b_{03})) + (b_{02}, b_{03}) \cdot L' \cr
&= (b_{02}^d - b_{03}^d, c_{11}(b_{02} - b_{11} b_{03})) + b_{03} \cdot L'
+ b_{02} L'', \cr
}
$$
where
$$
\eqalignno{
L' &= L'' +
c_{11} (b_{1i} - b_{1j}, b_{12}^d - 1)
+ \sum_{i=0}^{n-1} c_{11} \cdots c_{i1} \left(\D_{i+1} + B_i \right) \bigr), \cr
L'' &= (s - fb_{01}^d, b_{01}-b_{04}, b_{01}^d - b_{03}^d)
+ \D_0 + \D_1. \cr
}
$$
Note that  $L'$ is generated by all the generators of $\cap_{r \ge 1} p_r$
other than $b_{02}^d - b_{03}^d$.
Then the intersection of the last three components of $A$ is
$$
\eqalignno{
(c_{03}&,b_{01},b_{03},b_{02}^d) \cap (b_{02}, b_{03}) \cap \bigcap_{r = 1} ^n p_r \cr
&=
(b_{02}^d - b_{03}^d) + b_{03} \cdot L' + (c_{03},b_{01},b_{03},b_{02}^d) \cap
\left((c_{11}(b_{02} - b_{11} b_{03})) + b_{02} L''\right) \cr
&=
(b_{02}^d - b_{03}^d, b_{02}(b_{01}^d - b_{03}^d)) + b_{03} \cdot L' \cr
&\hskip1em
+(c_{03},b_{01},b_{03},b_{02}^d) \cap
\left((c_{11}(b_{02} - b_{11} b_{03}))
+ b_{02}((s - fb_{01}^d, b_{01}-b_{04})
+ \D_0 + \D_1) \right) \cr
&=
(b_{02}^d - b_{03}^d, b_{02}(b_{01}^d - b_{03}^d)) + b_{03} \cdot L' \cr
&\hskip1em
+(c_{03},b_{01}) \cdot
\left((c_{11}(b_{02} - b_{11} b_{03}))
+ b_{02}((s - fb_{01}^d, b_{01}-b_{04})
+ \D_0 + \D_1) \right) \cr
&\hskip1em
+(b_{03},b_{02}^d) \cap
\left((c_{11}(b_{02} - b_{11} b_{03}))
+ b_{02}((s - fb_{01}^d, b_{01}-b_{04})
+ \D_0 + \D_1) \right) \cr
&=
(b_{02}^d - b_{03}^d, b_{02}(b_{01}^d - b_{03}^d)) + b_{03} \cdot L' \cr
&\hskip1em
+(c_{03},b_{01}) \cdot
\left((c_{11}(b_{02} - b_{11} b_{03}))
+ b_{02}((s - fb_{01}^d, b_{01}-b_{04})
+ \D_0 + \D_1) \right) \cr
&\hskip1em
+(b_{03},b_{02}^d)
\cap
\left(c_{11}(b_{02} - b_{11} b_{03}),b_{02}\right)\cr
&\hskip2em \cap
(c_{11}(b_{02} - b_{11} b_{03}))
+ (s - fb_{01}^d, b_{01}-b_{04})
+ \D_0 + \D_1) \cr
&=
(b_{02}^d - b_{03}^d, b_{02}(b_{01}^d - b_{03}^d)) + b_{03} \cdot L' \cr
&\hskip1em
+(c_{03},b_{01}) \cdot
\left((c_{11}(b_{02} - b_{11} b_{03}))
+ b_{02}((s - fb_{01}^d, b_{01}-b_{04})
+ \D_0 + \D_1) \right) \cr
&\hskip1em
+(b_{03}c_{11}b_{11},b_{02}^d,b_{02}b_{03})
\cap
((c_{11}(b_{02} - b_{11} b_{03}),s - fb_{01}^d, b_{01}-b_{04})
+ \D_0 + \D_1) \cr
&=
(b_{02}^d - b_{03}^d, b_{02}(b_{01}^d - b_{03}^d)) + b_{03} \cdot L' \cr
&\hskip1em
+(c_{03},b_{01}) \cdot
\left((c_{11}(b_{02} - b_{11} b_{03}))
+ b_{02}((s - fb_{01}^d, b_{01}-b_{04})
+ \D_0 + \D_1) \right) \cr
&\hskip1em
+(b_{03}c_{11}b_{11},b_{02}^d,b_{02}b_{03})
\cdot ((s - fb_{01}^d, b_{01}-b_{04})+ \D_0+ \D_1) \cr
&\hskip1em
+(b_{03}c_{11}b_{11},b_{02}^d,b_{02}b_{03})
\cap
(c_{11}(b_{02} - b_{11} b_{03})) \cr
&=
(b_{02}^d - b_{03}^d, b_{02}(b_{01}^d - b_{03}^d)) + b_{03} \cdot L' \cr
&\hskip1em
+(c_{03},b_{01}) \cdot
\left((c_{11}(b_{02} - b_{11} b_{03}))
+ b_{02}((s - fb_{01}^d, b_{01}-b_{04})
+ \D_0 + \D_1) \right) \cr
&\hskip1em
+c_{11}(b_{02} - b_{11} b_{03})(b_{03},b_{02}^{d-1}) \cr
&=
(b_{02}^d - b_{03}^d, b_{02}(b_{01}^d - b_{03}^d),b_{02}^{d-1}c_{11}(b_{02} - b_{11} b_{03}))
+ b_{03} \bigcap_{r=1}^n p_r \cr
&\hskip1em
+(c_{03},b_{01}) \cdot
\left((c_{11}(b_{02} - b_{11} b_{03}))
+ b_{02}((s - fb_{01}^d, b_{01}-b_{04})
+ \D_0 + \D_1) \right). \cr
}
$$
Hence $A$, the intersection of all of its components,
equals
$$
\eqalignno{
A &= (c_{03},b_{01},c_{1i} b_{1i}^d,b_{02}^d) \cap (c_{03},b_{01},b_{03},b_{02}^d) \cap (b_{02}, b_{03})
\cap \bigcap_{r = 1} ^n p_r \cr
&=
(c_{03},b_{01}) \cdot
\left((c_{11}(b_{02} - b_{11} b_{03}))
+ b_{02}((s - fb_{01}^d, b_{01}-b_{04})
+ \D_0 + \D_1) \right) \cr
&\hskip1em
+(c_{03},b_{01},c_{1i} b_{1i}^d,b_{02}^d)
\cap
\bigl( (b_{02}^d - b_{03}^d, b_{02}(b_{01}^d - b_{03}^d),b_{02}^{d-1}c_{11}(b_{02} - b_{11} b_{03}))
+ b_{03} \bigcap_{r=1}^n p_r\bigr)  \cr
&=
(c_{03},b_{01}) \cdot
\left((c_{11}(b_{02} - b_{11} b_{03}))
+ b_{02}((s - fb_{01}^d, b_{01}-b_{04})
+ \D_0 + \D_1) \right) \cr
&\hskip1em
+c_{03}
\bigl( (b_{02}^d - b_{03}^d, b_{02}(b_{01}^d - b_{03}^d),b_{02}^{d-1}c_{11}(b_{02} - b_{11} b_{03}))
+ b_{03} \bigcap_{r=1}^n p_r\bigr)  \cr
&\hskip1em
+(b_{01},c_{1i} b_{1i}^d,b_{02}^d)
\cap
\bigl( (b_{02}^d - b_{03}^d, b_{02}(b_{01}^d - b_{03}^d),b_{02}^{d-1}c_{11}(b_{02} - b_{11} b_{03}))
+ b_{03} \bigcap_{r=1}^n p_r\bigr).  \cr
}
$$
Let $A'$ be the ideal in the last line.
The second intersectand ideal of $A'$ decomposes as
$$
\eqalignno{
(b_{02}^d &- b_{03}^d, b_{02}(b_{01}^d - b_{03}^d),b_{02}^{d-1}c_{11}(b_{02} - b_{11} b_{03}))
+ b_{03} \bigcap_{r=1}^n p_r  \cr
&=
\bigcap_{r=1}^n p_r 
\cap
\bigl((b_{02}^d,b_{03}^d, b_{02}(b_{01}^d - b_{03}^d),b_{02}^{d-1}c_{11}(b_{02} - b_{11} b_{03}))
+ b_{03} \bigcap_{r=1}^n p_r \bigr), \cr
%&=
%b_{03} \bigcap_{r=1}^n p_r
%+\bigl(b_{02}^d-b_{03}^d, b_{02}(b_{01}^d - b_{03}^d),
%b_{02}^{d-1}c_{11}(b_{02} - b_{11} b_{03}) \bigr)
%+ \bigcap_{r=1}^n p_r 
%\cap
%\bigl(b_{02}^d\bigr) \cr
%&=
%b_{03} \bigcap_{r=1}^n p_r
%+\bigl(b_{02}^d-b_{03}^d, b_{02}(b_{01}^d - b_{03}^d),
%b_{02}^{d-1}c_{11}(b_{02} - b_{11} b_{03}) \bigr)
% \cr
}
$$
so that
$$
\eqalignno{
A' &= \bigcap_{r=1}^n p_r
\cap (b_{01},c_{1i} b_{1i}^d,b_{02}^d)
\cap
\bigl((b_{02}^d,b_{03}^d, b_{02}b_{01}^d,b_{02}^{d-1}c_{11}b_{11} b_{03})
+ b_{03} \bigcap_{r=1}^n p_r \bigr)  \cr
&=\bigcap_{r=1}^n p_r
\cap
\left((b_{02}^d,b_{02}b_{01}^d)
+
(b_{01},c_{1i} b_{1i}^d,b_{02}^d)
\cap
\bigl((b_{03}^d, b_{02}^{d-1}c_{11}b_{11} b_{03})
+ b_{03} \bigcap_{r=1}^n p_r \bigr) \right) \cr
&=\bigcap_{r=1}^n p_r
\cap
\left((b_{02}^d,b_{02}b_{01}^d)
+
b_{03} \left((b_{01},c_{1i} b_{1i}^d,b_{02}^d)
\cap
\bigl((b_{03}^{d-1}, b_{02}^{d-1}c_{11}b_{11})
+ \bigcap_{r=1}^n p_r \bigr) \right) \right) \cr
&=\bigcap_{r=1}^n p_r
\cap
\left((b_{02}^d,b_{02}b_{01}^d,
b_{03}(c_{1i} b_{1i}^d-c_{11} b_{11}^d), b_{03}b_{01}^d) \right.\cr
&\hskip2em \left.
+
b_{03} \left((b_{01},c_{11} b_{11}^d)
\cap
\bigl((b_{03}^{d-1}, b_{02}^{d-1}c_{11}b_{11})
+ \bigcap_{r=1}^n p_r \bigr) \right) \right) \cr
&=\bigcap_{r=1}^n p_r
\cap
\left((b_{02}^d,b_{02}b_{01}^d,
b_{03}(c_{1i} b_{1i}^d-c_{11} b_{11}^d), b_{03}b_{01}^d) \right.\cr
&\hskip2em \left.
+
b_{03} \left((b_{01},c_{11} b_{11}^d)
\cap
\bigl((b_{03}^{d-1}, b_{02}^{d-1}c_{11}b_{11})
+ L''' \bigr) \right) \right), \cr
}
$$
where $L'''$ is generated by all the given generators of $\bigcap_{r \ge 1} p_r$
other than $b_{01}^d-b_{03}^d$:
$$
L''' = (s-fb_{01}^d, b_{01}-b_{04},b_{02}^d-b_{03}^d)
+ \D_0 + \C_1 F + \sum_{i=0}^{n-1} \C_1 \cdots \C_i (\D_{i+1} + B_i).
$$
With this,
$$
\eqalignno{
A' &=\bigcap_{r=1}^n p_r
\cap
\left((b_{02}^d,b_{02}b_{01}^d,
b_{03}(c_{1i} b_{1i}^d-c_{11} b_{11}^d), b_{03}b_{01}^d)\right.\cr
&\hskip2em
+
b_{03} b_{01}\bigl((b_{03}^{d-1}, b_{02}^{d-1}c_{11}b_{11})+ L''' \bigr) \left.
+ b_{03} \left((c_{11} b_{11}^d)
\cap
\bigl((b_{03}^{d-1}, b_{02}^{d-1}c_{11}b_{11})
+ L''' \bigr) \right) \right) \cr
&=\bigcap_{r=1}^n p_r
\cap
\left((b_{02}^d,b_{02}b_{01}^d,
b_{03}(c_{1i} b_{1i}^d-c_{11} b_{11}^d), b_{03}b_{01}^d)\right.\cr
&\hskip2em
+
b_{03} b_{01}\bigl((b_{03}^{d-1}, b_{02}^{d-1}c_{11}b_{11})+ L''' \bigr) \left.
+ b_{03}c_{11} b_{11}^d \left(
(b_{03}^{d-1}, b_{02}^{d-1})+ L''' : c_{11} \right)\right) \cr
&=
(b_{02}(b_{01}^d-b_{02}^d),
b_{03}(c_{1i} b_{1i}^d-c_{11} b_{11}^d),
b_{03}(b_{01}^d-b_{02}^d),
b_{01}(b_{03}^d-b_{02}^d)) \cr
&\hskip1em
+(b_{01}b_{02}^{d-1}c_{11}(b_{02}-b_{03}b_{11}))
+b_{03} b_{01}L''' \cr
&\hskip1em
+ (c_{11}b_{11}^d(b_{03}^d - b_{02}^d),
b_{02}^{d-1}b_{11}^{d-1}c_{11}(b_{02}-b_{03}b_{11}))
+ b_{03}c_{11} b_{11}^d \left(L''' : c_{11} \right)
+\bigcap_{r=1}^n p_r
\cap\left(b_{02}^d \right) \cr
&=
(b_{02}(b_{01}^d-b_{02}^d),
b_{03}(c_{1i} b_{1i}^d-c_{11} b_{11}^d),
b_{03}(b_{01}^d-b_{02}^d),
b_{01}(b_{03}^d-b_{02}^d)) \cr
&\hskip1em
+(b_{01}b_{02}^{d-1}c_{11}(b_{02}-b_{03}b_{11}))
+b_{03} b_{01}L''' \cr
&\hskip1em
+ (c_{11}b_{11}^d(b_{03}^d - b_{02}^d),
b_{02}^{d-1}b_{11}^{d-1}c_{11}(b_{02}-b_{03}b_{11}))
+ b_{03}c_{11} b_{11}^d \left(L''' : c_{11} \right)
+b_{02}^d \cdot \bigcap_{r=1}^n p_r. \cr
}
$$
Hence $A$ equals
$$
\eqalignno{
A &=
(c_{03},b_{01}) \cdot
\left((c_{11}(b_{02} - b_{11} b_{03}))
+ b_{02}((s - fb_{01}^d, b_{01}-b_{04})
+ \D_0 + \D_1) \right) \cr
&\hskip1em
+c_{03}
\bigl( (b_{02}^d - b_{03}^d, b_{02}(b_{01}^d - b_{03}^d),b_{02}^{d-1}c_{11}(b_{02} - b_{11} b_{03}))
+ b_{03} \bigcap_{r=1}^n p_r\bigr)  \cr
&\hskip1em
+(b_{02}(b_{01}^d-b_{02}^d),
b_{03}(c_{1i} b_{1i}^d-c_{11} b_{11}^d),
b_{03}(b_{01}^d-b_{02}^d),
b_{01}(b_{03}^d-b_{02}^d)) \cr
&\hskip1em
+(b_{01}b_{02}^{d-1}c_{11}(b_{02}-b_{03}b_{11}))
+b_{03} b_{01}L''' \cr
&\hskip1em
+ (c_{11}b_{11}^d(b_{03}^d - b_{02}^d),
b_{02}^{d-1}b_{11}^{d-1}c_{11}(b_{02}-b_{03}b_{11}))
+ b_{03}c_{11} b_{11}^d \left(L''' : c_{11} \right)
+b_{02}^d \cdot \bigcap_{r=1}^n p_r. \cr
}
$$
Thus finally,
$$
\eqalignno{
&\bigcap_{r = -4}^n p_r=
J + s \D_0 + f c_{02} A \cr
&=
J + s \D_0 + 
f c_{02}(c_{03},b_{01}) \cdot
\left(b_{02}((s - fb_{01}^d, b_{01}-b_{04})+ \D_0 + \D_1) \right) \cr
&\hskip1em
+f c_{02}c_{03}
\bigl( (b_{02}^d - b_{03}^d, b_{02}(b_{01}^d - b_{03}^d))
+ b_{03} \bigcap_{r=1}^n p_r\bigr)  \cr
&\hskip1em
+f c_{02}(b_{02}(b_{01}^d-b_{02}^d),
b_{03}(c_{1i} b_{1i}^d-c_{11} b_{11}^d),
b_{03}(b_{01}^d-b_{02}^d),
b_{01}(b_{03}^d-b_{02}^d),c_{11}b_{11}^d(b_{03}^d - b_{02}^d)) \cr
&\hskip1em
+f c_{02}b_{03} b_{01}L'''
+ f c_{02}b_{03}c_{11} b_{11}^d \left(L''' : c_{11} \right)
+f c_{02}b_{02}^d \cdot \bigcap_{r=1}^n p_r, \cr
}
$$
or in nicer form:

\thm
The intersection of all the minimal components of $J(n,d)$ equals
$$
\eqalignno{
\bigcap_{r = -4} ^n p_r
&= J + s\D_0
+ fc_{02}b_{02}(c_{03},b_{01})
\left(s - fb_{01}^d, b_{01}-b_{04}\right) \cr
&\hskip2em
+ fc_{02}b_{02}(c_{03},b_{01})(\D_0 + \D_1) \cr
&\hskip2em
+ fc_{02}c_{03}(b_{02}^d - b_{03}^d, b_{02}(b_{01}^d - b_{03}^d)) \cr
&\hskip2em
+ fc_{02}c_{03}b_{03}\cdot (E + \D_0 + \C_1 F +\sum_{i=0}^{n-1} c_{11} \cdots c_{i1} \left(\D_{i+1} + B_i \right)
) \cr
&\hskip2em
+f c_{02}(b_{02}(b_{01}^d-b_{02}^d),
b_{03}(c_{1i} b_{1i}^d-c_{11} b_{11}^d),
b_{03}(b_{01}^d-b_{02}^d)) \cr
&\hskip2em
+f c_{02}(b_{01}(b_{03}^d-b_{02}^d),c_{11}b_{11}^d(b_{03}^d - b_{02}^d)) \cr
&\hskip2em
+f c_{02}b_{03} b_{01}(
E'''
+ \D_0 + \C_1 F +\sum_{i=0}^{n-1} c_{11} \cdots c_{i1} \left(\D_{i+1} + B_i \right)
) \cr
&\hskip2em
+fc_{02}b_{03} b_{11}^d c_{11}(E'''
+ \D_0 + F +\sum_{i=0}^{n-1} c_{21} \cdots c_{i1} \left(\D_{i+1} + B_i \right)
) \cr
&\hskip2em
+f c_{02}b_{02}^d \cdot (E + \D_0 + \C_1 F +\sum_{i=0}^{n-1} c_{11} \cdots c_{i1} \left(\D_{i+1} + B_i \right)
), \cr
}
$$
where
$$
E''' = \left(s - fb_{01}^d, b_{01}-b_{04},b_{03}^d-b_{02}^d\right).
\eqed
$$
\endb

(With Macaulay2 I verified this theorem and intermediate computations
in the proof above for the case $n = 3$, $d = 2$.)

Set $c' = c_{11}c_{21} \cdots c_{n-2,1} (c_{n-1,1} - c_{n-1,4})$.
With the listed generators,
together with the generators $h_{rj}$ of $J$,
$$
\eqalignno{
s&c_{01} \cdots c_{n-2,1} (c_{n-1,1} - c_{n-1,4})
= c_{01} (s - fb_{01}^d) c'
+ (f c_{01} -s c_{02}) b_{01}^d c' \cr
&\hskip1em
+ c_{02}(s - fb_{02}^d) b_{01}^d c'
+ b_{01}^d fc_{02}b_{02}^dc', \cr
}
$$
%
%%+ fc_{02}b_{02}^d b_{01}^d c'
%+ fc_{02}c_{11}(b_{02}-b_{11}b_{03})
%(b_{02}^{d-1} + b_{02}^{d-2} b_{11} b_{03}+ \cdots + b_{03}^{d-1} b_{11}^{d-1}) b_{01}^d c'' \cr
%&\hskip1em
%+ fc_{02}c_{11}b_{03}^d b_{11}^d b_{01}^d c'' \cr
%&= h_{01} c'
%+ h_{13} b_{01}^d c'
%+ h_{02} b_{01}^d c'
%+ h_{17}(b_{02}^{d-1} + b_{02}^{d-2} b_{11} b_{03}+ \cdots + b_{03}^{d-1} b_{11}^{d-1}) b_{01}^d c'' \cr
%&\hskip 1em
%+ \left(fc_{02}b_{03}b_{11}^d c_{11}\cdots c_{n-2,1} (c_{n-1,1} - c_{n-1,4})\right)
%b_{01}^d c'' b_{03}^{d-1}
%, \cr
%}
%$$
and the degrees of the coefficients $c'$, $c'$, $ b_{01}^d c'$ and $b_{01}^d fc_{02}b_{02}^d$
of the generators
$h_{01}, h_{13}, h_{02}$ and $fc_{02}b_{02}^d c'$, respectively,
of the intersection of the minimal components,
are not doubly exponential in $n$.
This proves:

\thm
The doubly exponential ideal membership problem of the Mayr-Meyer ideals $J(n,d)$ and $J_l(n,d)$
for the element
$sc_{01} \cdots c_{n-2,1} (c_{n-1,1} - c_{n-1,4})$
is not due to the minimal components,
but to some embedded prime ideal.
\qed
\endb

%________________________________________________________
\vskip 4ex

%\vfill\eject
\bigskip
\leftline{\bf References}
\bigskip

\bgroup
\font\eightrm=cmr8 \def\rm{\fam0\eightrm}
\font\eightit=cmti8 \def\it{\fam\itfam\eightit}
\font\eightbf=cmbx8 \def\bf{\fam\bffam\eightbf}
\font\eighttt=cmtt8 \def\tt{\fam\ttfam\eighttt}
\rm
\baselineskip=9.9pt
\parindent=3.6em

\item{[BS]}
D.\ Bayer and M.\ Stillman,
On the complexity of computing syzygies,
{\it J.\ Symbolic Comput.}, {\bf 6} (1988), 135-147.

\item{[D]}
M.\ Demazure,
Le th\'eor\`eme de complexit\'e de Mayr et Meyer,
{\it G\'eom\'etrie alg\'ebrique et applications, I
(La R\'abida, 1984)}, 35-58,
Travaux en Cours, 22, {\it Hermann, Paris}, 1987.

\item{[EHV]}
D.\ Eisenbud, C.\ Huneke and W.\ Vasconcelos,
Direct methods for primary decomposition,
{\it Invent.\ math.}, {\bf 110} (1992), 207-235.

\item{[ES]}
D.\ Eisenbud and B.\ Sturmfels,
Binomial ideals,
{\it Duke Math.\ J.}, {\bf 84} (1996), 1-45.

\item{[GTZ]}
P.\ Gianni, B.\ Trager and G.\ Zacharias,
Gr\"obner bases and primary decompositions of polynomial ideals,
{\it J.\ Symbolic Comput.}, {\bf 6} (1988), 149-167.

\item{[GS]}
D.\ Grayson and M.\ Stillman,
Macaulay2. 1996.
A system for computation in algebraic geometry and commutative algebra,
available via anonymous {\tt ftp} from {\tt math.uiuc.edu}.

\item{[GPS]}
G.-M.\ Greuel, G.\ Pfister and H.\ Sch\"onemann,
Singular. 1995.
A system for computation in algebraic geometry and singularity theory.
Available via anonymous {\tt ftp} from {\tt helios.mathematik.uni-kl.de}.

\item{[H]}
G.\ Herrmann,
Die Frage der endlich vielen Schritte in der Theorie der Polynomideale,
{\it Math.\ Ann.}, {\bf 95} (1926), 736-788.

\item{[K]}
J.\ Koh,
Ideals generated by quadrics exhibiting double exponential degrees,
{\it J.\ Algebra}, {\bf 200} (1998), 225-245.

\item{[MM]}
E.\ Mayr and A.\ Meyer,
The complexity of the word problems for commutative semigroups
and polynomial ideals,
{\it Adv.\ Math.}, {\bf 46} (1982), 305-329.

\item{[SY]}
T.\ Shimoyama and K.\ Yokoyama,
Localization and primary decomposition of polynomial ideals,
{\it J.\ of Symbolic Comput.}, {\bf 22} (1996), 247-277.

\item{[S1]}
I.\ Swanson,
The first Mayr-Meyer ideal,
preprint, 2001.

\item{[S2]}
I.\ Swanson,
On the embedded primes of the Mayr-Meyer ideal,
preprint, 2002.

\item{[S3]}
I.\ Swanson,
A new family of ideals with the doubly exponential ideal membership property,
preprint, 2002.

\egroup

\vskip 4ex
\noindent
{\sl
New Mexico State University - Department of Mathematical Sciences,
Las Cruces, New Mexico 88003-8001, USA.
E-mail: {\tt iswanson@nmsu.edu}.
}

\end